\documentclass[11pt]{article}
\usepackage{amscd}
\usepackage{amsfonts}
\usepackage{amsmath}
\usepackage{amssymb}
\usepackage{amsthm}
\usepackage{bbm}
\usepackage{CJK}
\usepackage{fancyhdr}
\usepackage{graphicx}
\usepackage{hyperref}
\usepackage{indentfirst}
\usepackage{latexsym}
\usepackage{mathrsfs}
\usepackage{xypic}
\usepackage[compress]{cite}
\newtheorem{theorem}{Theorem}[section]
\newtheorem{lemma}[theorem]{Lemma}
\newtheorem{definition}[theorem]{Definition}
\newtheorem{proposition}[theorem]{Proposition}
\newtheorem{example}[theorem]{Example}

\newtheorem{corollary}[theorem]{Corollary}
\newtheorem{remark}[theorem]{Remark}
\usepackage[top=1in,bottom=1in,left=1.25in,right=1.25in]{geometry}
\def\<{\langle}
\def\>{\rangle}
\def\a{\alpha}
\def\b{\beta}

\def\c{\cdot}

\def\d{\delta}

\def\g{\gamma}

\date{}
\begin{document}
\renewcommand{\baselinestretch}{1.2}
\renewcommand{\arraystretch}{1.0}
\title{\bf On split regular Hom-Leibniz-Rinehart algebras}
\author{{\bf Shuangjian Guo$^{1}$, Xiaohui Zhang$^{2}$,  Shengxiang Wang$^{3}$\footnote
        { Corresponding author(Shengxiang Wang):~~wangsx-math@163.com} }\\
{\small 1. School of Mathematics and Statistics, Guizhou University of Finance and Economics} \\
{\small  Guiyang  550025, P. R. of China} \\
{\small 2.  School of Mathematical Sciences, Qufu Normal University}\\
{\small Qufu  273165, P. R. of China}\\
{\small 3.~ School of Mathematics and Finance, Chuzhou University}\\
 {\small   Chuzhou 239000,  P. R. of China}}
 \maketitle
\begin{center}
\begin{minipage}{13.cm}

{\bf \begin{center} ABSTRACT \end{center}}
In this paper,  we introduce the notion of the Hom-Leibniz-Rinehart algebra as an algebraic analogue of Hom-Leibniz algebroid,  and prove that such an arbitrary split regular Hom-Leibniz-Rinehart algebra $L$ is of the form  $L=U+\sum_{\g}I_\g$ with $U$ a subspace of a maximal abelian subalgebra $H$ and any $I_{\g}$, a well described ideal of $L$, satisfying $[I_\g, I_\d]= 0$ if
$[\g]\neq [\d]$.  In the sequel,  we  develop techniques of connections of roots and weights for split Hom-Leibniz-Rinehart algebras respectively. Finally, we study the structures of tight split regular Hom-Leibniz-Rinehart algebras.


 \end{minipage}
 \end{center}
 \normalsize\vskip1cm

\section*{INTRODUCTION}
\def\theequation{0. \arabic{equation}}
\setcounter{equation} {0}
The notion of the Lie-Rinehart algebra plays an important role in many branches of mathematics. The idea of
this notion goes back to the work of Jacobson to study certain field extensions.
It was also appeared in some
different names in several areas which includes differential geometry and differential Galois theory. In
\cite{Mackenzie05}, Mackenzie provided a list of 14 different terms mentioned for this notion. Huebschmann viewed Lie-Rinehart algebras as an algebraic
counterpart of Lie algebroids defined over smooth manifolds.  His work on several aspects of this algebra
has been developed systematically through a series of articles namely \cite{Huebschmann90,Huebschmann98, Huebschmann99, Huebschmann04} .

The notion of Hom-Lie algebras was first introduced by Hartwig,
Larsson and Silvestrov in \cite{Hartwig}, who developed an approach to deformations 
of the Witt and  Virasoro algebras basing on $\sigma$-deformations.
In fact, Hom-Lie algebras include Lie algebras as a subclass, but the deformation  of Lie algebras twisted by a homomorphism.

 Mandal and  Mishra defined modules over a Hom-Lie-Rinehart algebra and studied a cohomology with coefficients
in a left module. They presented the notion of extensions of Hom-Lie-Rinehart
algebras and deduced a characterisation of low dimensional cohomology spaces in terms of the
group of automorphisms of certain abelian extensions and the equivalence classes of those
abelian extensions in the category of Hom-Lie-Rinehart algebras in \cite{Mandal2017}. The concept of a Hom-Lie-Rinehart algebra  has a geometric analogue which is nowadays
called a Hom-Lie algebroid in \cite{Cai17} and \cite{Laurent-Gengoux}. See also \cite{Castiglioni18,Mandal2018, Mandal17, Mandal18,  zhang18} for other works on Hom-Lie-Rinehart algebras.

The class of the split algebras is specially related to addition quantum numbers, graded contractions and deformations. For instance, for a physical system
which displays a symmetry, it is interesting to know the detailed structure of the split
decomposition, since its roots can be seen as certain eigenvalues which are the additive
quantum numbers characterizing the state of such system. Determining the structure of
split algebras will become more and more meaningful in the area of research in mathematical physics. Recently, in \cite{Aragon2015}-\cite{Cao18}, the structure of different classes of split algebras have been determined by the techniques of connections of roots.  Recently,   we studied  the structures of split regular Hom-Lie Rinehart algebras in \cite{Wang19}.     The purpose of this paper is to  consider the structure of split regular Hom-Leibniz-Rinehart algebras by the
techniques of connections of roots based on some work in \cite{Cao18} and  \cite{Wang19} .

This paper is organized
as follows.  In Section 2, we introduce the notion of the Hom-Leibniz-Rinehart algebra and prove that such an arbitrary split regular Hom-Leibniz-Rinehart algebras $L$ is of the form  $L=U+\sum_{\g}I_\g$ with $U$ a subspace of a maximal abelian subalgebra $H$ and any $I_{\g}$, a well described ideal of $L$, satisfying $[I_\g, I_\d]= 0$ if
$[\g]\neq [\d]$.  In Section 3 and Section 4,  we  develop techniques of connections of roots and weights for split Hom-Leibniz-Rinehart algebras respectively. In the last section, we study the structures of tight split regular Hom-Leibniz-Rinehart algebras.

\section{Preliminaries}
\def\theequation{\arabic{section}.\arabic{equation}}
\setcounter{equation} {0}

Let $R$ denote a commutative ring with unity,   $\mathbb{Z}$ the set of all integers  and $\mathbb{N}$ the set of all nonnegative integers,  all algebraic systems are considered of arbitrary dimension and over an arbitrary base field  $\mathbb{K}$.     And we recall some basic definitions and results related to our paper from \cite{Makhlouf10} and \cite{Mandal2017}.

\begin{definition}
Given an associative commutative algebra $A$, an $A$-module $M$  and an algebra
endomorphism $\phi: A \rightarrow  A$, we call an $R$-linear map $\delta: A \rightarrow M$ a $\phi$-derivation of $A$ into  $M$ if it satisfies
the required identity:
\begin{eqnarray*}
\delta(ab)=\phi(a)\delta(b)+\phi(b)\delta(a), ~~~\mbox{for any  $a, b\in A$}.
\end{eqnarray*}
\end{definition}
Let us denote by $Der_{\phi}(A)$ the set of $\phi$-derivations of $A$ into itself.

\begin{definition}
A Hom-Leibniz algebra $L$ is an algebra $L$, endowed with a bilinear
product
\begin{eqnarray*}
[\c,\c]: L\times L\rightarrow L,
\end{eqnarray*}
and a homomorphism $\psi: L\rightarrow L$
\begin{eqnarray*}
&&[\psi(x), [y, z]]=[[x,y],\psi(z)]+[\psi(y), [x,z]]~~\mbox{(Hom-Leibniz~identity)}
\end{eqnarray*}
holds for any $x, y, z\in L$.
\end{definition}
If $\psi$ is furthermore an algebra automorphism, that is,  a linear bijective on such that $\psi([x,y])=[\psi(x),\psi(y)]$ for any $x,y\in L$, then $L$ is called a regular Hom-Leibniz algebra.

\begin{definition}
A Hom-Lie-Rinehart algebra over $(A, \phi)$ is a tuple $(A, L,  [\c, \c], \phi, \psi, \rho)$, where $A$ is an
associative commutative algebra, $L$ is an $A$-module, $[\c, \c] :  L \times L \rightarrow L$ is a skew symmetric bilinear map,
$\phi: A\rightarrow A$ is an algebra homomorphism, $\psi: L\rightarrow L$ is a linear map satisfying $\psi([x, y])=[\psi(x), \psi(y)]$,
and the $R$-map $\rho: L\rightarrow Der_{\phi}(A)$  such that following conditions hold.\\
(1) The triple $(L,  [\c, \c], \psi)$ is a Hom-Lie algebra.\\
(2) $\psi(a\c x)=\phi(a)\c \psi(x)$ for all $a\in A, x \in L$.\\
(3) $(\rho, \phi)$ is a representation of $(L,  [\c, \c], \psi)$ on $A$.\\
(4) $\rho(a\c x)=\phi(a)\c \rho(x)$ for all $a\in A, x \in L$.\\
(5) $[x, a\c y] = \phi(a)\c [x, y] + \rho(x)(a)\psi(y)$ for all $a\in A, x, y\in L$.
\end{definition}
A Hom-Lie-Rinehart algebra $(A, L,  [\c, \c], \phi, \psi, \rho)$ is said to be regular if the map $\phi: A\rightarrow A$ is
an algebra automorphism and $\psi: L\rightarrow L$ is a bijective map.

\section{Decomposition}
\def\theequation{\arabic{section}. \arabic{equation}}
\setcounter{equation} {0}

In this section,  we introduce the notion of the Hom-Leibniz-Rinehart algebra as an algebraic analogue of
Hom-Leibniz algebroid. In a sequel, we introduce the class of split algebras in the framework of Hom-Leibniz-Rinehart algebras.

\begin{definition}
A Hom-Leibniz-Rinehart algebra over $(A, \phi)$ is a tuple $(A, L,  [\c, \c], \phi, \psi, \rho)$, where $A$ is an
associative commutative algebra, $L$ is an $A$-module, $[\c, \c] : L\times L \rightarrow L$ is a skew symmetric linear map,
$\phi: A\rightarrow A$ is an algebra homomorphism, $\psi: L\rightarrow L$ is a linear map satisfying $\psi([x, y])=[\psi(x), \psi(y)]$,
and the $R$-map $\rho: L\rightarrow Der_{\phi}(A)$  satisfying the  following conditions.\\
(1) The triple $(L,  [\c, \c], \psi)$ is a Hom-Leibniz algebra.\\
(2) $\psi(a\c x)=\phi(a)\c \psi(x)$ for all $a\in A, x \in L$.\\
(3) $(\rho, \phi)$ is a representation of $(L,  [\c, \c], \psi)$ on $A$.\\
(4) $\rho(a\c x)=\phi(a)\c \rho(x)$ for all $a\in A, x\in L$.\\
(5) $[x, a\c y] = \phi(a)\c [x, y] + \rho(x)(a)\psi(y)$ for all $a\in A, x, y\in L$.
\end{definition}

We denote it by $(L,A)$ or just by $L$ if there is not any possible confusion.
A Hom-Leibniz-Rinehart algebra $(A, L,  [\c, \c], \phi, \psi, \rho)$ is said to be regular if the map $\phi: A\rightarrow A$ is
an algebra automorphism and $\psi: L\rightarrow L$ is a bijective map.

\begin{example}
A Leibniz-Rinehart algebra $L$ over $A$ with the linear map $[\c, \c] : L\times L  \rightarrow L$ and the  $R$-map $\rho: L\rightarrow Der(A)$ is a Hom-Leibniz-Rinehart algebra
$(A, L,  [\c, \c], \phi, \psi, \rho)$, where $\psi=Id_L, \phi=Id_A$  and $\rho: L\times L\rightarrow Der_{\phi}(A)=Der(A)$.

\end{example}
\begin{example}
A Hom-Leibniz algebra $(L,  [\c, \c],  \psi)$ structure over an $R$-module $L$ gives the Hom-Leibniz-Rinehart algebra $(A, L,  [\c, \c], \phi, \psi, \rho)$ with $A=R$, the algebra morphism $\phi=id_{R}$ and the trivial action of $L$ on $R$.

If we consider a Leibniz-Rinehart algebra $L$ over $A$ along with an endomorphism
\begin{eqnarray*}
(\phi, \psi):(A, L)\rightarrow (A, L)
\end{eqnarray*}
in the category of Leibniz-Rinehart algebras,  then we get a Hom-Leibniz-Rinehart algebra $(A, L,  [\c, \c]_\psi, \phi, \psi, \rho_{\phi})$ as follows:\\
(1) $[x, y]_{\psi}=\psi([x, y])$ for any $x, y\in L$.\\
(2) $\rho_{\phi}(x)(a)=\phi(\rho(x)(a))$ for all $a\in A, x\in L$.
\end{example}

\begin{definition}
A Hom-Leibniz algebroid is a tuple $(\xi,  [\c, \c], \phi, \psi, \rho)$, where $\xi: A\rightarrow M$ is a vector bundle over a smooth manifold $M$, $\phi: M\rightarrow M$ is a smooth map,
 $[-, -] : \Gamma(A)\times \Gamma(A) \rightarrow \Gamma(A)$ is a bilinear map,
the map $\rho: \phi^{!}A\rightarrow \phi^{!}TM$ is called the anchor and $\psi:\Gamma(A) \rightarrow \Gamma(A)$ is a linear map  such that following conditions are satisfied .\\
(1) The triplet $(\Gamma(A),  [-, -], \psi)$ is a Hom-Lie algebra.\\
(2) $\psi(fX)=\phi^{\ast}(F)\psi(X)$ for all $X\in \Gamma(A), f \in C^{\infty}(M)$.\\
(3) $(\rho, \phi^{\ast})$ is a representation of $(\Gamma(A),  [-, -], \psi)$ on $C^{\infty}(M)$.\\
(5) $[X, fY] = \phi^{\ast}(f) [X, Y] + \rho(X)[f]\psi(Y)$ for all $X, Y\in \Gamma(A), f\in C^{\infty}(M)$.
\end{definition}

\begin{example}
A Hom-Leibniz algebroid provides a Hom-Leibniz-Rinehart algebra $(C^{\infty}(M), \Gamma(A),$ \\ $[\c, \c], \phi^{\ast}, \psi, \rho)$, where $\Gamma(A)$ is the space  a sections of the underline vector bundle $A\rightarrow M$ and $\phi^{\ast}: C^{\infty}(M)\rightarrow C^{\infty}(M)$ is canonically defined by the smooth map $\phi: M\rightarrow M$.
\end{example}

\begin{example}
Let $(A, L,  [\c, \c]_L, \phi, \psi_L, \rho_L)$  and $(A, M,  [\c, \c]_M, \phi, \psi_M, \rho_M)$ be two Hom-Leibniz-Rinehart algebras over $(A, \phi)$. We consider
\begin{eqnarray*}
L\times_{Der_{\phi}A} M=\{(l, m)\in L\times M: \rho_L(l)=\rho_M(m)\}.
\end{eqnarray*}
Then $(A, L\times_{Der_{\phi}A} M, [\c, \c], \phi, \psi, \widetilde{\rho})$ is a  Hom-Leibniz-Rinehart algebra, where \\
(1)  The linear bracket $[\c, \c]$  is given by
\begin{eqnarray*}
[(l_1,m_1),(l_2,m_2)]:=([l_1, l_2], [m_1, m_2]),
\end{eqnarray*}
for any $l_1,l_2 \in L$ and $m_1,m_2\in M$.\\
(2) The map $\psi: L\times_{Der_{\phi}A} M\rightarrow L\times_{Der_{\phi}A} M$ is given by
\begin{eqnarray*}
\psi(l, m):=(\psi_L(l), \psi_M(m)),
\end{eqnarray*}
for any $l\in L$ and $m\in M$.\\
(3) The action of $L\times_{Der_{\phi}A} M$ on $A$ is given by
\begin{eqnarray*}
\widetilde{\rho}(l, m)(a):=\rho_L(l)(a)=\rho_M(m)(a),
\end{eqnarray*}
for any $l\in L, m\in M$ and $a\in A$.
\end{example}
Next we define homomorphisms of Hom-Leibniz-Rinehart algebras.
\begin{definition}
Let $(A, L,  [\c, \c]_L, \phi, \psi_L, \rho_L)$ and  $(B, L',  [\c, \c]_{L'}, \psi, \psi_{L'}, \rho_{L'})$ be  two Hom-Leibniz-Rinehart algebras, then a Hom-Leibniz-Rinehart algebra homomorphism is defined as a pair of maps $(g, f)$, where the map $g: A \rightarrow B$ is an $R$-algebra homomorphism and $f: L\rightarrow L'$ is an $R$-linear map such that
following identities hold:\\
(1) $f(a\c x) =g(a)\c f(x)$, for all $x\in L$ and $a\in A$. \\
(2) $f([x, y]_L)=[f(x), f(y)]_{L'}$, for all $x, y\in L$.\\
(3) $f(\psi_L(x))=\psi_{L'}(f(x))$, for all $x\in L$.\\
(4) $g(\phi(a))=\psi(g(a))$, for all $a\in A$.\\
(5) $g(\rho_L(x)(a))=\rho_{L'}(f(x))(g(a))$, for all $x\in L$ and $a\in A$.
\end{definition}

\begin{definition}
 A subalgebra $(S,A)$ of $(L,A)$ is called a \emph{Hom-Leibniz subalgebra},  if $(S,A)$ satisfies $AS\subset S$ such that $S$ acts on $A$ via the composition $S\hookrightarrow L\rightarrow Der_{\phi} (A)$.
A Hom-Leibniz subalgebra $(I,A)$ of $(L, A)$ is called an \emph{ideal},  if $I$ is a Hom-Leibniz ideal of $L$ such that $\rho(I)(A)L\subseteq I.$
\end{definition}

The ideal $J$ generated by
\begin{eqnarray*}
\{[x, y]+[y, x]: x, y\in L\}
\end{eqnarray*}
plays an important role in mathematics since it determines the non-super Lie character of $L$. From Hom-Leibniz identity, it is straightforward to check that this ideal satisfies
\begin{eqnarray}
[L, J]=0.
\end{eqnarray}

Let us introduce the class of split algebras in the framework of Hom-Leibniz algebras from \cite{Cao18}. Denote by $H$ a maximal abelian subalgebra of a Hom-Leibniz algebra $L$.  For a linear functional
\begin{eqnarray*}
\g:H\rightarrow \mathbb{K},
\end{eqnarray*}
we define the root space of $L$ associated to $\g$ as the subspace
\begin{eqnarray*}
L_{\g}:=\{v_{\a}\in L:[h,\psi(v_{\g})]=\a(h)\psi(v_{\g}), \mbox{for any $h\in H$}\}.
\end{eqnarray*}
The elements $\g:H\rightarrow \mathbb{K}$ satisfying $L_{\g}\neq 0$ are called roots of $L$ with respect to $H$ and we denote $\Gamma:=\{\g\in H^{\ast}\setminus\{0\}:L_{\g}\neq 0\}$. We call that $L$ is a split regular Hom-Leibniz algebra with respect to $H$ if
\begin{eqnarray*}
L=H\oplus \bigoplus_{\g\in \Gamma}L_{\g}.
\end{eqnarray*}
We also say that $\Gamma$ is the root system of $L$.

\begin{definition}
 A \emph{ split regular Hom-Leibniz-Rinehart algebra} (with respect to a MASA $H$ of the regular Hom-Leibniz algebra $L$, here MASA means maximal abelian subalgebra)
  is a regular Hom-Leibniz-Rinehart algebra $(L,A)$ in which the Hom-Leibniz algebra $L$ contains a splitting Cartan
subalgebra $H$ and the algebra $A $ is a weight module (with respect to $H$) in the sense that
$A$ can be written as the direct sum
$A=A_0\oplus (\bigoplus_{\alpha\in \Lambda}A_{\alpha})$ with $\phi(A_\alpha)\subset A_\alpha$,
where
\begin{eqnarray*}
A_{\alpha}:=\{a_{\alpha}\in A|\rho(h)(a_{\alpha})=\alpha(h)\phi(a_{\alpha}),~\forall h\in H\},
\end{eqnarray*}
for a linear functional $\alpha\in H^{\ast}$ and $\Lambda:=\{\alpha\in H^{\ast}\backslash\{0\}: A_\a\neq 0\}$
 denotes the weights system of $A$. The linear subspace $ A_{\alpha}$ , for $\alpha\in \Lambda$, is called the\emph{ weight space of
$A$ associate to $\alpha$}, the element $\alpha\in \Lambda\cup \{0\}$ are called \emph{weights} of $A$.
\end{definition}

\begin{example}
Let $(L,A)$ be a Leibniz-Rinehart algebra, where  $L=H\oplus (\bigoplus_{\a\in \Gamma}L_{\g}), ~A=A_0\oplus (\bigoplus_{\alpha\in \Lambda}A_{\alpha})$ , $\psi: L\rightarrow L,  \phi: A\rightarrow A$ are two  automorphisms such that $\psi(H)=H, \phi(A_0)=A_0$ and  $\phi(A_\alpha)\subset A_\alpha$. By Example  2.3, we know that $(L, A)$ is a regular Hom-Leibniz-Rinehart algebra. Then we have
\begin{eqnarray*}
L=H\oplus (\bigoplus_{\a\in \Gamma}\mathfrak{L}_{\a\psi^{-1}}),~~~   A=A_0\oplus (\bigoplus_{\alpha\in \Lambda}A_{\alpha\phi^{-1}}),
\end{eqnarray*}
which make the regular Hom-Leibniz-Rinehart algebra $(L,  A)$ being the roots system $\Gamma'=\{\a\psi^{-1}:\a\in \Gamma\}$ and weights system $\Lambda':=\{\alpha\phi^{-1}:\a\in \Lambda\}$.
\end{example}

The following  lemma is  analogous  to the results  of \cite{Wang19}.

\begin{lemma}
For any $ \gamma,\xi\in\Gamma\cup \{0\}$  and $ \alpha,\beta\in\Lambda\cup \{0\}$,
 the following assertions hold.

(1) $L_0=H.$

(2) $\psi(L_{\g})=L_{\g\psi^{-1}}$ and $\psi^{-1}(L_{\g})=L_{\g\psi}$.

(3) If $[L_{\gamma},L_{\xi}]\neq 0$, then
 $ \gamma\psi^{-1}+\xi\psi^{-1}\in\Gamma\cup \{0\}$ and $[L_{\gamma},L_{\xi}]\subset L_{\gamma\psi^{-1}+\xi\psi^{-1}}$.

(4) If $A_{\alpha}A_{\beta}\neq 0$, then
 $\alpha+\beta\in\Lambda\cup \{0\}$ and $A_{\alpha}A_{\beta}\subset A_{\alpha+\beta}$.

(5) If $A_{\alpha}L_{\gamma}\neq 0$, then
 $\alpha+\gamma\in\Gamma\cup \{0\}$ and $A_{\alpha}L_{\gamma}\subset L_{\alpha+\gamma}$.

(6) If $\rho(L_{\gamma})A_{\alpha}\neq 0$, then
$\alpha+\gamma\in\Lambda\cup \{0\}$ and $\rho(L_{\gamma})A_{\alpha}\subset A_{\alpha+\gamma}$.

\end{lemma}

\section{Connections of roots}
\def\theequation{\arabic{section}. \arabic{equation}}
\setcounter{equation} {0}

In what follows, $L$ denotes a split regular Hom-Leibniz-Rinehart algebra and
\begin{eqnarray*}
L=H\oplus (\bigoplus_{\g\in \Gamma}L_{\g}), \ \ \  A=A_0\oplus (\bigoplus_{\alpha\in \Lambda}A_{\alpha}).
\end{eqnarray*}
Given a linear functional $\gamma: H\rightarrow  \mathbb{K}$, we denote by $-\gamma: H\rightarrow  \mathbb{K}$ the element
in $H^{\ast}$ defined by $(-\gamma)(h):=-\gamma(h)$ for all $h\in  H.$ We also denote $-\Gamma:=\{-\gamma: \gamma\in \Gamma\}$.
In a similar way we can define $-\Lambda:=\{-\alpha: \alpha\in\Lambda\}$.
Finally, we denote $\pm\Gamma:=\Gamma\cup -\Gamma $ and $\pm\Lambda:=\Lambda\cup -\Lambda.$

\begin{definition}
Let $\gamma,\xi\in \Gamma$,
 we say that $\gamma$  is  connected to $\xi$ if

 $\bullet$ Either $\xi=\epsilon\gamma\psi^{z}$ for some $z\in \mathbb{Z}$ and $\epsilon\in \{1,-1\}$.

  $\bullet$ Either there exists a family $\{\zeta_1,\zeta_2,\ldots,\zeta_n\}\subset\pm\Lambda\cup\pm\Gamma$,
     with $n\geq 2$, such that

(1) $\zeta_1\in\{\gamma\psi^{k}|k\in \mathbb{Z}\}$.

(2) $\zeta_1\psi^{-1}+\zeta_2\psi^{-1}\in\pm\Gamma$,

    $~~~~~\zeta_1\psi^{-2}+\zeta_2\psi^{-2}+\zeta_3\psi^{-1}\in\pm\Gamma$,

    $~~~~~\zeta_1\psi^{-3}+\zeta_2\psi^{-3}+\zeta_3\psi^{-2}+\zeta_4\psi^{-1}\in\pm\Gamma$,

    $~~~~~~~~~~\cdots\cdots\cdots$

    $~~~~~\zeta_1\psi^{-i}+\zeta_2\psi^{-i}+\zeta_3\psi^{-i+1}+\cdots+\zeta_{i+1}\psi^{-1}\in\pm\Gamma$,

    $~~~~~~~~~~\cdots\cdots\cdots$

    $~~~~~\zeta_1\psi^{-n+2}+\zeta_2\psi^{-n+2}+\zeta_3\psi^{-n+3}+\cdots+\zeta_{n-1}\psi^{-1}\in\pm\Gamma$.

(3) $\zeta_1\psi^{-n+1}+\zeta_2\psi^{-n+1}+\zeta_3\psi^{-n+2}+\cdots+\zeta_{n}\psi^{-1}\in\{\pm\xi\psi^{-m}|m\in \mathbb{Z}\}$.

We will also say that $\{\zeta_1,\zeta_2,\ldots,\zeta_n\}$ is a \emph{connection from $\gamma$ to $\xi$}.
\end{definition}

The proof of the next result is analogous to the one of \cite{Cao18}.

\begin{proposition}
The relation $\sim$ in $\Gamma$ is an equivalence relation, where $\gamma\sim\xi$ if and only if $\gamma$  is  connected to $\xi$.
\end{proposition}

  By  Proposition 3.2 we can consider the quotient set
  \begin{eqnarray*}
\Gamma/\sim=\{[\g]:\g\in \Gamma\},
  \end{eqnarray*}
with $[\g]$ being the set of nonzero roots which are connected to $\g$.
Our next goal is to associate an ideal $I_{[\g]}$ to $[\g]$. Fix $[\g]\in \Gamma/\sim$, we start by defining
\begin{eqnarray*}
L_{0,[\g]}:=(\bigoplus_{\xi\in [\g], -\xi\in \Lambda}A_{-\xi}L_{\xi})+(\bigoplus_{\xi\in [\g]}[L_{-\xi}, L_{\xi}]).
\end{eqnarray*}
Now we define
\begin{eqnarray*}
L_{[\g]}:=\bigoplus_{\xi\in [\gamma]}L_{\xi}.
\end{eqnarray*}
Finally, we denote by $I_{[\gamma]}$ the direct sum of the two subspaces above:
\begin{eqnarray*}
I_{[\g]}:=L_{0,[\g]}\oplus L_{[\g]}.
\end{eqnarray*}

\begin{proposition} For any $[\g]\in \Lambda/\sim$, the following assertions hold.

(1) $[I_{[\g]},I_{[\g]}]\subset I_{[\g]}$.

(2) $\psi(I_{[\g]})=I_{[\g]}$.

(3) $AI_{[\g]}\subset I_{[\g]}$.

(4) $\rho(I_{[\g]})(A)L\subset I_{[\g]}$.

(5) For any $[\g]\neq [\delta]$, we have $[I_{[\g]}, I_{[\delta]}]=0$.
\end{proposition}
{\bf Proof.} (1) First we check that $[I_{[\g]},I_{[\g]}]\subset I_{[\g]}$, we can write
\begin{eqnarray}
 [I_{[\g]},I_{[\g]}]&=&[L_{0,[\g]}\oplus L_{[\g]},L_{0,[\g]}\oplus L_{[\g]}]\nonumber\\
 &\subset&[L_{0,[\g]}, L_{[\g]}]+[L_{[\g]}, L_{0,[\g]}]+[L_{[\g]}, L_{[\g]}].
\end{eqnarray}
Given $\delta\in [\g]$, we have $[L_{0,[\g]}, L_{\delta}]\subset L_{\delta}\subset L_{[\g]}$.  By a similar argument, we get $[L_{\delta}, L_{0,[\g]}]\subset L_{[\g]}$.

Next we consider  $[L_{[\g]},L_{[\g]}]$. If we take $\delta,\eta\in [\g]$  such that $[L_{\delta}, L_{\eta}]\neq 0$, then $[L_{\delta}, L_{\eta}]\subset L_{\delta+\eta}$.
If $\delta\psi^{-1}+\eta\psi^{-1}=0$, we get $[L_{\delta}, L_{-\delta}]\subset L_{0,[\g]}$ . Suppose that $\delta\psi^{-1}+\eta\psi^{-1}\in \Gamma$. We infer that $\{\delta,\eta\}$ is a connection from $\delta$ to $\delta\psi^{-1}+\eta\psi^{-1}$. The transitivity of $\sim$ now gives that $\delta\psi^{-1}+\eta\psi^{-1}\in [\g]$ and so $[L_{\delta}, L_{\eta}]\subset L_{[\g]}$. Hence
\begin{eqnarray}
[L_{[\g]},L_{[\g]}]\in I_{[\g]}.
\end{eqnarray}
 From (3.1) and (3.2), we get $[I_{[\g]},I_{[\g]}]\subset I_{[\g]}$.

(2) It is easy to check that $\psi(I_{[\g]})=I_{[\g]}$.

(3) and (4) similar to \cite{Wang19}.

(5) We will study the expression $[I_{[\g]},I_{[\delta]}]$. Notice that
\begin{eqnarray}
 [I_{[\g]},I_{[\delta]}]&=&[L_{0,[\g]}\oplus L_{[\g]},L_{0,[\delta]}\oplus L_{[\delta]}]\nonumber\\
 &\subset&[L_{0,[\g]}, L_{[\delta]}]+[L_{[\g]}, L_{0,[\delta]}]+[L_{[\g]}, L_{[\delta]}].
\end{eqnarray}

First we  consider $[L_{[\gamma]}, L_{[\delta]}]$ and suppose that there exist $\g_1\in [\g], \d_1\in [\d]$ such that $[L_{\g_1}, L_{\d_1}]\neq 0$. As necessarily $\g_1\psi^{-1}\neq-\d_1\psi^{-1}$, then $\g_1\psi^{-1}+\d_1\psi^{-1}\in \Gamma$. So $\{\g_1,\d_1, -\g_1\psi^{-1}\}$ is a connection between $\g_1$ and $\d_1$. By the transitivity of the connection relation we see $\g\in [\d]$, a contradicition. Hence $[L_{\g_1}, L_{\d_1}]=0$, and so
\begin{eqnarray}
[L_{[\gamma]}, L_{[\delta]}]=0.
\end{eqnarray}
By the definition of $L_{0,[\gamma]}$, we have
$$[L_{0,[\gamma]},L_{[\delta]}]=[(\sum_{\gamma_1\in [\gamma],-\gamma_1\in\Lambda}A_{-\gamma_1}L_{\gamma_1})
+(\sum_{\gamma_{1}\in [\gamma]}[L_{-\gamma_1},L_{\gamma_1}]),L_{[\delta]}].$$

Suppose that there exist $\g_1\in [\g]$ and $\d_1\in [\d]$ such that
\begin{eqnarray*}
[L_{\d_1}, [L_{\g_1}, L_{-\g_1}]]= 0.
\end{eqnarray*}
Suppose that $
[L_{\d_1}, [L_{\g_1}, L_{-\g_1}]]\neq 0,
$
then  Hom-Leibniz~ identity gives
\begin{eqnarray*}
0&\neq& [ \psi\psi^{-1}(L_{\d_1}), [L_{\g_1}, L_{-\g_1}]]\\
&\subset&[[\psi^{-1}(L_{\d_1}), L_{\g_1}], \psi(L_{-\g_1})]+[[\psi^{-1}(L_{\d_1}), L_{-\g_1}], \psi(L_{\g_1})].
\end{eqnarray*}
Hence
\begin{eqnarray*}
[\psi^{-1}(L_{\d_1}), L_{\g_1}]+[\psi^{-1}(L_{\d_1}), L_{-\g_1}]\neq 0,
\end{eqnarray*}
which contradicts (3.4). Therefore, $[L_{\d_1}, [L_{\g_1}, L_{-\g_1}]]=0$.

For the expression  $[A_{-\gamma_1}L_{\gamma_1},L_{[\delta]}]$, suppose there exists $\delta_1\in [\delta]$ such that $[A_{-\gamma_1}L_{\gamma_1},L_{\delta_1}]\neq 0$.
 By Definition 2.1, we have
$$[A_{-\gamma_1}L_{\gamma_1},L_{\delta_1}]
=[L_{\delta_1},A_{-\gamma_1}L_{\gamma_1}]
\subset \phi(A_{-\gamma_1})[L_{\delta_1},L_{\gamma_1}]+\rho(L_{\delta_1})(A_{-\gamma_1})L_{\delta_{1}\psi^{-1}}.$$
By the discussion above, we get $[L_{\delta_1},L_{\gamma_1}]=0$.
Since $[A_{-\gamma_1}L_{\gamma_1},L_{[\delta]}]\neq 0$,
it follows that $0 \neq\rho(L_{\delta_1})(A_{-\gamma_1})L_{\delta_{1}\psi^{-1}}\subset A_{\delta_{1}-\gamma_{1}}L_{\delta_{1}\psi^{-1}}$.
Thus $A_{\delta_{1}-\gamma_{1}}\neq 0$ and $\delta_{1}-\gamma_{1}\in \Lambda\cup\{0\}$. $\delta_{1}\sim\gamma_{1}$, a contradiction.
So  $[A_{-\gamma_1}L_{\gamma_1},L_{[\delta]}]=0$.  Therefore, we have $[L_{0,[\gamma]},L_{[\delta]}]=0$. In a similar way we can prove  $[L_{[\g]}, L_{0,[\delta]}]=0$, we conclude $[I_{[\g]},I_{[\delta]}]=0$.  \hfill $\square$

\begin{definition}
 A Hom-Leibniz-Rinehart algebra $(L,A)$ is \emph{simple} if $[L,L]\neq 0,AA\neq 0,AL\neq 0$
  and its only ideals are $\{0\}, J, L $ and   the kernel of $\rho$.
\end{definition}

\begin{theorem}
The following assertions hold.

(1) For any $[\g]\in \Gamma/\sim$, the linear space $I_{[\g]}=L_{0,[\g]}+L_{[\g]}$ of $L$ associated to $[\g]$ is an ideal of $L$.

(2) If $L$ is simple, then there exists a connection from $\g$ to $\delta$ for any $\g,\delta\in \Gamma$ and \begin{eqnarray*}
H=(\sum_{\gamma\in \Gamma,-\gamma\in\Lambda}A_{-\gamma}L_{\gamma})+(\sum_{\gamma\in \Gamma}[L_{-\gamma},L_{\gamma}]).
\end{eqnarray*}
.
\end{theorem}
{\bf Proof.} (1) Since $[I_{[\g]},H]+[H, I_{[\g]}]\subset I_{[\g]}$, by Proposition 3.3, we have
\begin{eqnarray*}
[I_{[\g]}, L]=[I_{[\g]}, H\oplus(\bigoplus_{\xi\in [\g]} L_{\xi})\oplus (\bigoplus_{\delta\notin [\g]} L_{\delta})]\subset I_{[\g]}
\end{eqnarray*}
and
\begin{eqnarray*}
[L, I_{[\g]},]=[ H\oplus(\bigoplus_{\xi\in [\g]} L_{\xi})\oplus (\bigoplus_{\delta\notin [\g]} L_{\delta}), I_{[\g]}]\subset I_{[\g]}.
\end{eqnarray*}
Furthermore,
\begin{eqnarray*}
[I_{[\g]}, L]+[L, I_{[\g]}]
=[I_{[\g]}, H\oplus(\bigoplus_{\xi\in [\g]} L_{\xi})\oplus (\bigoplus_{\delta\notin [\g]} L_{\delta})]
+[H\oplus(\bigoplus_{\xi\in [\g]} L_{\xi})\oplus (\bigoplus_{\delta\notin [\g]} L_{\delta}), I_{[\g]}]\subset I_{[\g]}.
\end{eqnarray*}
As we also have $\psi(I_{[\g]})=I_{[\g]}$. So we show that $I_{[\gamma]}$ is a Hom-Leibniz ideal of $L$, we also have that $I_{[\gamma]}$ is an $A$-module, the we conclude $I_{[\gamma]}$ is an ideal of $L$.

(2) The simplicity of $L$ implies $I_{[\g]}\in \{J, L$, ker$\rho\}$. If some $\g \in \Gamma$ is such that $I_{[\g]}=L$, then $[\g]=\Gamma$.  Otherwise, if $I_{[\g]}=J$ for any $\a \in \Gamma$, then $[\g]=[\delta]$ for any $\g,\delta\in \Gamma$, and so $[\g]=\Gamma$.   Otherwise, if $I_{[\g]}= ker \rho$ for any $\g \in \Gamma$, then $[\g]=[\delta]$ for any $\g,\delta\in \Gamma$, and so $[\g]=\Gamma$. Thus $H=(\sum_{\gamma\in \Gamma,-\gamma\in\Lambda}A_{-\gamma}L_{\gamma})+(\sum_{\gamma\in \Gamma}[L_{-\gamma},L_{\gamma}])$.
\hfill $\square$

\begin{theorem} We have
\begin{eqnarray*}
L=U+\sum_{[\g]\in \Lambda/\sim}I_{[\g]},
\end{eqnarray*}
where $U$ is a linear complement in $H$ of $(\sum_{\gamma\in \Gamma,-\gamma\in\Lambda}A_{-\gamma}L_{\gamma})+(\sum_{\gamma\in \Gamma}[L_{-\gamma},L_{\gamma}])$ and any $I_{[\g]}$ is one of the ideals of $L$ described in Theorem 3.5, satisfying $[I_{[\g]},I_{[\delta]}]=0$ if $[\g]\neq[\delta]$.
\end{theorem}
{\bf Proof.}  $I_{[\g]}$ is well defined and is an ideal of $L$ and it is clear that
 \begin{eqnarray*}
L=H\oplus\sum_{[\g]\in \Gamma} L_{[\g]}=U+\sum_{[\g]\in \Gamma/\sim}I_{[\g]}.
 \end{eqnarray*}
 Finally, Proposition 3.3 gives us $[I_{[\g]},I_{[\delta]}]=0$ if $[\g]\neq[\delta]$.\hfill $\square$

\begin{definition}
The annihilator of a Hom-Leibniz-Rinehart algebra $L$ is the set
\begin{eqnarray*}
Z(L):=\{v\in L:[v,L]+[L, v]=0\mbox{~and~} \rho(v)=0\}.
\end{eqnarray*}
\end{definition}

 \begin{corollary}  If $Z(L)=0$ and $H=(\sum_{\gamma\in \Gamma,-\gamma\in\Lambda}A_{-\gamma}L_{\gamma})+(\sum_{\gamma\in \Gamma}[L_{-\gamma},L_{\gamma}])$. Then $L$ is the direct sum of the ideals given in Theorem 3.6,
 \begin{eqnarray*}
L=\bigoplus_{[\g]\in \Gamma/\sim}I_{[\g]},
\end{eqnarray*}
Furthermore, $[I_{[\g]},I_{[\delta]}]=0$ if $[\g]\neq[\delta]$.\hfill $\square$.
\end{corollary}
{\bf Proof.} Since $H=(\sum_{\gamma\in \Gamma,-\gamma\in\Lambda}A_{-\gamma}L_{\gamma})+(\sum_{\gamma\in \Gamma}[L_{-\gamma},L_{\gamma}])$, it follows that  $L=\sum_{[\g]\in \Gamma/\sim}I_{[\g]}$. To verify the direct character of the sum, take some $v\in I_{[\g]}\cap(\sum_{[\delta]\in\Gamma/\sim,[\delta]\neq[\g]}I_{[\delta]})$. Since $v\in I_{[\g]}$, the fact $[I_{[\g]},I_{[\delta]}]=0$ when  $[\g]\neq[\delta]$ gives us
\begin{eqnarray*}
[v, \sum_{[\delta]\in\Gamma/\sim,[\delta]\neq[\g]}I_{[\delta]}]+[\sum_{[\delta]\in\Gamma/\sim,[\delta]\neq[\g]}I_{[\delta]},v]=0.
\end{eqnarray*}
  In a similar way,  $v\in \sum_{[\delta]\in\Gamma/\sim,[\delta]\neq[\g]}I_{[\delta]}$ implies $[v,I_{[\g]}]+[I_{[\g]}, v]=0$.  It is easy to obtain that $ \rho(v)=0$.   That is $v\in Z(L)$ and so $v=0$.              \hfill $\square$

\section{Decompositions of  $A$}
\def\theequation{\arabic{section}. \arabic{equation}}
\setcounter{equation} {0}

  We will discuss the weight  spaces and decompositions of $A$ similar to \cite{Wang19} and omit the proof.

\begin{definition}
Let $\alpha,\beta\in \Lambda$
 we say that $\alpha$ and $\beta$ are \emph{connected} if

 $\bullet$ either $\beta=\varepsilon\alpha$ for some   $\varepsilon\in \{1,-1\}$;

  $\bullet$ either there exists a family $\{\sigma_1,\sigma_2,\ldots,\sigma_n\}\subset\pm\Lambda\cup\pm\Gamma$,
     with $n\geq 2$, such that

(1) $\sigma_1=\alpha$.

(2) $\sigma_1+\sigma_2\in\pm\Lambda\cup\pm\Gamma$,

    $~~~~~\sigma_1+\sigma_2+\sigma_3\in\pm\Lambda\cup\pm\Gamma$,

    $~~~~~~~~~~\cdots\cdots\cdots$

    $~~~~~\sigma_1+\sigma_2\cdots+\sigma_{n-1}\in\pm\Lambda\cup\pm\Gamma$,

(3) $\sigma_1+\sigma_2\cdots+\sigma_{n}\in\{\beta,-\beta\}$.

We will also say that $\{\sigma_1,\sigma_2,\ldots,\sigma_n\}$ is a \emph{connection} from $\alpha$ to $\beta$.
\end{definition}

\begin{proposition}
The relation $\approx$ in $\Lambda$ is an equivalence relation, where $\alpha\approx\beta$ if and only if
 $\alpha$  is  connected to $\beta$.
\end{proposition}

By Proposition 4.2, we can consider the quotient set
$$\Lambda/\approx:=\{[\alpha]|\alpha\in \Lambda\},$$
where $[\alpha]$ denotes the set of nonzero weights of $A$ which are connected to $\alpha$.
In the following we will  associate an  adequate  ideal $\mathcal{A}_{[\alpha]}$ to any $[\alpha]$.
 For a fixed $\alpha\in \Lambda$, we define
\begin{eqnarray*}
A_{0,[\alpha]}:=\bigl(\sum_{\beta\in [\alpha],-\beta\in\Lambda}\rho(L_{-\beta})(A_{\beta})
                    +(\sum_{\beta\in [\alpha]}A_{-\beta},A_{\beta})\bigr)\subset A_0,~~
 A_{[\alpha]}:=\bigoplus_{\beta\in [\alpha]}A_{\beta}.
 \end{eqnarray*}
Then we denote by $\mathcal{A}_{[\alpha]}$ the direct sum of the two subspaces above,
 \begin{eqnarray*}
 \mathcal{A}_{[\alpha]}:=A_{0,[\alpha]}\oplus A_{[\alpha]}.
 \end{eqnarray*}

\begin{proposition}
For any $\alpha,\beta\in \Lambda$, the following assertions hold.

(1) $\mathcal{A}_{[\alpha]}\mathcal{A}_{[\alpha]}\subset \mathcal{A}_{[\alpha]}$.

 (2) If $[\alpha]\neq[\beta]$, then $\mathcal{A}_{[\alpha]}\mathcal{A}_{[\beta]}=0$.
\end{proposition}

\begin{theorem}
Let $A$ be a commutative and associative algebra associated to a Hom-Leibniz-Rinehart algebra $L$.
Then the following assertions hold.

(1) For any $[\alpha]\in\Lambda/\approx$, the linear space
$\mathcal{A}_{[\alpha]}=A_{0,[\alpha]}\oplus A_{[\alpha]}$
of $A$ associated to $[\alpha]$ is an ideal of $A$.

(2) If $A$ is simple, then all weights of $\Lambda$ are connected. Furthermore,
$$A_{0}=\sum_{-\alpha\in \Gamma,\alpha\in\Lambda}\rho(L_{-\alpha})(A_{\alpha})
                    +(\sum_{\alpha\in \Lambda}A_{-\alpha}A_{\alpha}).$$
\end{theorem}

\begin{theorem}
Let $A$ be a commutative and associative algebra associated to a Hom-Leibniz-Rinehart algebra $L$.
Then
\begin{eqnarray*}
A=V+\sum_{[\alpha]\in\Lambda/\approx}\mathcal{A}_{[\alpha]},
\end{eqnarray*}
where $V$ is a linear complement in $A_0$ of
$\sum_{-\alpha\in \Gamma,\alpha\in\Lambda}\rho(L_{-\alpha})(A_{\alpha})
                    +(\sum_{\alpha\in \Lambda}A_{-\alpha}A_{\alpha})$
and any $A_{[\alpha]}$ is one of the ideals of $A$ described in Theorem 4.4-(1),
satisfying $\mathcal{A}_{[\alpha]}\mathcal{A}_{[\beta]}=0$, whenever $[\alpha]\neq[\beta]$.
\end{theorem}

Let us denote by $Z(A)$ the center of $A$, that is,
$Z(A):=\{a\in L| aA=0\}.$

\begin{corollary}
Let $(L,A)$ be a  Hom-Leibniz-Rinehart algebra.
If $Z(A)=0$ and
$$A_{0}=\sum_{-\alpha\in \Gamma,\alpha\in\Lambda}\rho(L_{-\alpha})(A_{\alpha})
                    +(\sum_{\alpha\in \Lambda}A_{-\alpha}A_{\alpha}),$$
then $A$ is the direct sum of the ideals given in Theorem 4.5, that is,
    \begin{eqnarray*}
A=\sum_{[\alpha]\in\Lambda/\approx}\mathcal{A}_{[\alpha]},
\end{eqnarray*}
  satisfying $\mathcal{A}_{[\alpha]}\mathcal{A}_{[\beta]}=0$, whenever $[\alpha]\neq[\beta]$.
\end{corollary}

\section{The simple components}
\def\theequation{\arabic{section}. \arabic{equation}}
\setcounter{equation} {0}

In this section we focus on the simplicity of split regular Hom-Leibniz-Rinehart algebra $(L,A)$ by centering our attention in those of maximal length.
From now on we always assume that  $\Lambda$ is symmetric in the sense that $\Lambda=-\Lambda$.
 \medskip

\begin{lemma}
 Let $(L,A)$ be a split regular Hom-Leibniz-Rinehart and $I$ an ideal of $L$.
Then $I=(I\cap H)\oplus(I\cap \bigoplus_{\gamma\in\Gamma}L_\gamma)$.
\end{lemma}
{\bf Proof.}
 Since $(L,A)$ is split, we get $L=H\oplus(\bigoplus_{\gamma\in\Gamma}L_\gamma)$.
 By the assumption that $I$ is an ideal of $L$, it is clear that $I$ is a submodule of $L$. Since a submodule of a weight module is again a
weight module.
 Thus $I$ is a  weight module and therefore
$I=(I\cap H)\oplus(I\cap \bigoplus_{\gamma\in\Gamma}L_\gamma)$.$\hfill \Box$

\begin{lemma}
 Let $(L,A)$ be a split regular Hom-Leibniz-Rinehart algebra with $Z(L)=0$ and $I$ an ideal of $L$.
If $I \subset H$, then $I=\{0\}$.
\end{lemma}
{\bf Proof.}
Since  $I \subset H$, $[I,H]+[H,I]\subset [H,H]=0$.
It follows that $[I,L]+[L, I]=[I,\bigoplus_{\gamma\in\Gamma}L_\gamma]+[\bigoplus_{\gamma\in\Gamma}L_\gamma, I]\subset H\cap (\bigoplus_{\gamma\in\Gamma}L_\gamma)=0$.
So  $I \subset Z(L)=0$.
$\hfill \Box$
\medskip

Observe that if $L$ is of  maximal length, then we have
\begin{eqnarray}
I=(I\cap H)\oplus (\bigoplus_{\g\in\Gamma^{I}} L_{\g})),
\end{eqnarray}
where $\Gamma^{I}=\{\g\in \Gamma: I\cap L_{\g}\neq 0\}$.

In particular, in case $I=J$, we get
\begin{eqnarray}
J=(J\cap H)\oplus (\bigoplus_{\g\in \Gamma^{J}}  L_{\g})
\end{eqnarray}
with $\Gamma^{J}=\{\g\in \Gamma: J\cap L_{\g}\neq 0\}=\{\g\in \Gamma: 0\neq L_{\g}\subset J\}$.

From here, we can write
\begin{eqnarray}
\Gamma=\Gamma^{J}\cup \Gamma^{\neg J},
\end{eqnarray}
where $
    \Gamma^{\neg J}=\{\g\in \Gamma: L_{\g}\neq 0 ~~\mbox{and}~~ J\cap L_{\g}= 0\}.
      $
      Therefore, we can write
      \begin{eqnarray}
L=H\oplus (\bigoplus_{\g\in \Gamma^{J}} L_{\g})\oplus (\bigoplus_{\delta\in \Gamma^{\neg J}}  L_{\delta}).
      \end{eqnarray}

Let us introduce the notion of root-multiplicativity in the framework of split
regular Hom-Leibniz-Rinehart algebras of maximal length, in a similar way to the ones for
split regular Hom-Lie Rinehart algebras in \cite{Wang19}.
\begin{definition}
A split regular Hom-Leibniz-Rinehart algebra
$(L,A)$ is called \emph{root-multiplicative} if  the following conditions hold:

(1) Given $\gamma,\delta\in\Gamma^{\neg J}$ such that  $\gamma\psi^{-1}+\delta\psi^{-1}\in\Gamma$, then $[L_{\gamma},L_{\delta}]\neq 0.$

(2) Given $\gamma \in\Gamma^{J}, \delta\in\Gamma^{\neg J}$ such that  $\gamma\psi^{-1}+\delta\psi^{-1}\in \Gamma^{J}$, then $[L_{\gamma},L_{\delta}]\neq 0.$

(3) Given $\a\in \Lambda, \gamma \in\Gamma$ such that $\alpha+\gamma\in\Gamma$, then $A_\alpha L_{\gamma}\neq 0.$

(4) If $\alpha+\beta\in\Lambda$, then $A_\alpha A_\beta\neq 0.$
\end{definition}
\begin{definition}
A split regular Hom-Leibniz-Rinehart $(L,A)$ is called \emph{of maximal length} if
dim $L_{\gamma}$=dim $A_\alpha=1$  for any $\gamma\in\Gamma$ and $\alpha\in\Lambda$.
\end{definition}

\begin{proposition}
Suppose   $ H=(\sum_{\gamma\in\Gamma^{\neg J},-\gamma\in\Lambda}A_{-\gamma}L_{\gamma})+(\sum_{\gamma\in \Gamma^{\neg J}}[L_{-\gamma},L_{\gamma}])$, $Z_{Lie}(L)=0$ and  root-multiplicative. If  $\Gamma^{\neg J}$ has all of its roots $\neg J$-connected,  then any ideal $I$ of $L$ such that $I\nsubseteq H\oplus J$, then $I=L$.
\end{proposition}
{\bf Proof.} By (5.1) and (5.3), we can write
\begin{eqnarray*}
I=(I\cap H)\oplus (\bigoplus_{\g\in \Gamma^{\neg J,I}} L_{\g}),
\end{eqnarray*}
where $\Gamma^{\neg J,I}=\Gamma^{\neg J}\cap \Gamma^{I}$ and $\Gamma^{J,I}=\Gamma^{J}\cap \Gamma^{I}$. Since $I\nsubseteq H\oplus J$, there exists $\g_0\in \Gamma^{\neg J}$ such that
\begin{eqnarray}
0\neq L_{\g_0} \subset I.
\end{eqnarray}
By Lemma 2.11, $\psi(L_{\g_0})=L_{\g_0\psi^{-1}}$. Equation(5.5) gives us $\psi(L_{\g_0}) \subset \psi(I)=I$. So $L_{\g_0\psi^{-1}}\subset I$. Similarly we get
\begin{eqnarray}
L_{\g_0\psi^{-n}}\subset I,  ~~\mbox{for}~ n\in \mathbb{N}.
\end{eqnarray}
For any $\delta\in \Gamma^{\neg J}$, $\b\notin \pm \g_0\psi^{-n}$, for  $n\in \mathbb{N}$, the fact that $\g_0$ and $\g$ are $\neg J$-connected gives us a $\neg J$-connection $\{\g_1,\g_2,...,\g_n\}\subset \Gamma^{\neg J}$ from $\g_0$ to $\delta$ such that\\
$\g_1=\g_0\in\Gamma^{\neg J},  \g_k\in \Gamma^{\neg J}$, for  $k=2,...,n$, \\
 $\g_1\psi^{-1}+\g_2\psi^{-1}\in \Gamma^{\Upsilon}$,\\
 $\cdots\cdots\cdots\cdots$\\
  $\g_1\psi^{-n+1}+\g_2\psi^{-n+1}+\g_3\psi^{-n+2}+\cdot\cdot\cdot+\g_{n-1}\psi^{-2}+\g_{n}\psi^{-1}\in \Gamma^{\Upsilon}$,\\
 $\g_1\psi^{-n+1}+\g_2\psi^{-n+1}+\g_3\psi^{-n+2}+\cdot\cdot\cdot+\g_{i}\psi^{-n+i-1}+\cdot\cdot\cdot+\g_{n-1}\psi^{-2}+\g_{n}\psi^{-1} \in \{\pm \b\phi^{-m}:m\in \mathbb{N}\}.$

 Taking into account $\g_1=\g_0\in\Gamma^{\neg J}$,  if $\g_2\in \Lambda$  (respectively $\g_2\in \Gamma^{\neg J}$),   the root-multiplicativity
and maximal length of $(L, A)$ allow us to assert
 \begin{eqnarray*}
0\neq A_{\g_1}L_{\g_2}=L_{\g_1+\g_2} ~(\mbox{respectively}, 0\neq [L_{\g_1}, L_{\g_2}]=L_{\g_1\psi^{-1}+\g_2\psi^{-1}}).
 \end{eqnarray*}
By (5.6), we have
 \begin{eqnarray*}
0\neq L_{\g_1\psi^{-1}+\g_2\psi^{-1}}\subset I.
 \end{eqnarray*}
 We can discuss in a similar way from $\g_1\psi^{-1}+\g_2\psi^{-1}\in \Gamma^{\neg J}, \g_3\in \Lambda \cup \Gamma^{\neg J}$ and $\g_1\psi^{-2}+\g_2\psi^{-2}+\g_3\psi^{-1}\in \Gamma^{\neg J}$ to get
  \begin{eqnarray*}
0\neq [L_{\g_1\psi^{-1}+\g_2\psi^{-1}}, L_{\g_3}]=L_{\g_1\psi^{-2}+\g_2\psi^{-2}+\g_3\psi^{-1}}.
 \end{eqnarray*}
Thus we have
 \begin{eqnarray*}
0\neq  L_{\g_1\psi^{-2}+\g_2\psi^{-2}+\g_3\psi^{-1}}\subset I.
 \end{eqnarray*}
 Following this process with the $\neg J$-connection $\{\g_1,\g_2,...,\g_n\}$ we obtain that
  \begin{eqnarray*}
0\neq L_{\g_1\psi^{-n+1}+\g_2\psi^{-n+1}+\g_3\psi^{-n+2}+...+\g_n\psi^{-1}}\subset I.
 \end{eqnarray*}
It follows that  either
 \begin{eqnarray}
L_{\delta\psi^{-m}}\subset I ~~\mbox{or}~~L_{-\delta\psi^{-m}}\subset I
 \end{eqnarray}
 for any $\delta\in \Gamma^{\neg J}, m\in \mathbb{N}$. Moreover, we have
 \begin{eqnarray}
\delta\psi^{-m}\in \Gamma^{\neg  J}.
 \end{eqnarray}
Since $ H=(\sum_{\gamma\in\Gamma^{\neg J},-\gamma\in\Lambda}A_{-\gamma}L_{\gamma})+(\sum_{\gamma\in \Gamma^{\neg J}}[L_{-\gamma},L_{\gamma}])$, by (5.7) and (5.8), we get
\begin{eqnarray}
H\subset I.
\end{eqnarray}
Now, for any $\Upsilon\in \{J, \neg J\}$, given any $\delta\in \Gamma^{J}$, the facts $\delta\neq0, H\subset I$ and the maximal length of $(L, A)$ show that
\begin{eqnarray*}
L_{\delta}=[L_{\delta\psi}, H]\subset I.
\end{eqnarray*}
The decomposition of $L$ in (5.4) finally gives us $H=I$. \hfill $\square$

Another interesting notion related to a split Hom-Leibniz-Rinehart algebra  of maximal length $(L, A)$ is Lie annihilator.
Write $L=H\oplus (\bigoplus_{\g\in \Gamma^{\neg J}} L_{\g})\oplus (\bigoplus_{\delta\in \Gamma^{ J}}  L_{\delta})$.

\begin{definition}
The Lie-annihilator of a split Hom-Leibniz-Rinehart algebra of maximal length $(L, A)$ is the set
\begin{eqnarray*}
Z_{Lie}(L):=\{v\in L:[v,H\oplus (\bigoplus_{\g\in \Gamma^{\neg J}} L_{\g})]+[H\oplus (\bigoplus_{\g\in \Gamma^{\neg J}} L_{\g}), v]=0\mbox{~and~} \rho(v)=0\}.
\end{eqnarray*}
\end{definition}
Observe that $Z(L)\subset Z_{Lie}(L)$.

In the following, we will discuss the relation between the decompositions of $L$ and $A$ of a Hom-Leibniz-Rinehart algebra $(L,A)$.

\begin{definition}
A split regular Hom-Leiniz-Rinehart algebra $(L,A)$ is tight if $Z_{Lie}(L)=0,Z(A)=0,AA=A,AL=L$ and
\begin{eqnarray*}
   H=(\sum_{\gamma\in \Gamma^{\neg J}, -\g\in \Lambda}A_{-\gamma}L_{\gamma})+(\sum_{\gamma\in \Gamma^{\neg J}}[L_{-\gamma},L_{\gamma}]),
~A_{0}=(\sum_{-\alpha\in \Gamma^{\neg J},\alpha\in\Lambda}\rho(L_{-\alpha})(A_{\alpha}))
                    +(\sum_{\alpha\in \Lambda}A_{-\alpha}A_{\alpha}).
\end{eqnarray*}
\end{definition}

\begin{remark}
Let $(L,A)$ be a tight split regular Hom-Leibniz-Rinehart algebra, then
    \begin{eqnarray*}
L=\sum_{[\gamma]\in\Gamma^{\neg J}/\sim}I_{[\gamma]},~A=\sum_{[\alpha]\in\Lambda/\approx}\mathcal{A}_{[\alpha]},
\end{eqnarray*}
with any $I_{[\gamma]}$ an ideal of $L$ verifying $[I_{[\gamma]},I_{[\delta]}]=0$ if $ [\gamma]\neq[\delta] $
and any $\mathcal{A}_{[\alpha]}$ an ideal of $A$
satisfying  $\mathcal{A}_{[\alpha]}\mathcal{A}_{[\beta]}=0$ if $ [\alpha]\neq[\beta].$
\end{remark}

\begin{proposition}
Let $(L,A)$ be a tight split regular Hom-Leibniz-Rinehart algebra, then for any $[\gamma]\in\Gamma^{\neg J}/\sim$
there exists a unique $[\alpha]\in\Lambda/\approx$ such that $\mathcal{A}_{[\alpha]}I_{[\gamma]}=0$.
\end{proposition}
{\bf Proof.}  Similar to Proposition 4.2 in \cite{Cao18}.
$\hfill \Box$

\begin{theorem}
Let $(L,A)$ be a tight split regular Hom-Leibniz-Rinehart algebra, then
    \begin{eqnarray*}
L=\sum_{i\in\Gamma^{\neg J}/ I}L_{i},~A=\sum_{j\in K}A_{j},
\end{eqnarray*}
with any $L_i$ a nonzero ideal of $L$ and any $A_j$ a nonzero ideal of $A$.
Furthermore, for any $i\in I$ there exists a unique $\tilde{j}\in K$ such that $A_{\tilde{j}}L_i=0$.
 \end{theorem}

\begin{theorem}
  Let $(L, A)$ be a tight split regular Hom-Leibniz-Rinehart algebra  of maximal length and root  multiplicative. If $\Gamma^{J}$,
$\Gamma^{\neg J}$ are symmetric and $\Gamma^{\neg J}$ has all of its roots $\neg J$-connected, then any ideal $I$ of $L$ such that $I\subseteq  J$ satisfies  either $I=J$ or $ J = I \oplus I'$ with $I'$ an ideal of $L$.
\end{theorem}
{\bf Proof.} By (5.1), we can write
\begin{eqnarray}
I=(I\cap H)\oplus (\bigoplus_{\g\in\Gamma^{J I}} L_{\g})),
\end{eqnarray}
and with $\Gamma^{J I}\subset \Gamma^{J}$. For any $\g\notin \Gamma^{J}$, we have
\begin{eqnarray*}
[J\cap H,  L_{\g}]+[L_{\g}, J\cap H]\subset L_{\g}\subset J.
\end{eqnarray*}
Hence, in case $[J\cap H,  L_{\g}]+[L_{\g}, J\cap H]\neq 0$ we have $\g\in \Gamma^{J}$, a contradiction. Hence  $[J\cap H,  L_{\g}]+[L_{\g}, J\cap H]= 0$, and so
\begin{eqnarray}
J\cap H \subset Z_{Lie}(L).
\end{eqnarray}
Taking into account $I\cap H\subset J\cap H=0$, we also write
\begin{eqnarray*}
I=\bigoplus_{\delta\in \Gamma^{J I}}  L_{\delta},
\end{eqnarray*}
with $\Gamma^{J I}\subset \Gamma^{J}$. Hence, we can take some $\delta_0\in \Gamma^{I}$ such that
\begin{eqnarray*}
0\neq L_{\delta_0}\subset I.
\end{eqnarray*}
Now, we can argue with the root-multiplicativity and the maximal length of $L$ as in Proposition 5.5 to
conclude that given any $\delta \in \Gamma^{J}$, there exists a $\neg J$-connection
$\{\delta_1,\delta_2, ..., \delta_r\}$ from $\delta_0$ to $\delta$ such that
\begin{eqnarray*}
0\neq [[...[L_{\delta_0}, L_{\delta_2}],...], L_{\delta_r}]\in L_{\delta\psi^{-m}}, ~\mbox{for}~m\in \mathbb{N}
\end{eqnarray*}
and so
\begin{eqnarray}
L_{\epsilon\delta\psi^{-m}}\subset I, ~~\mbox{for~some}~\epsilon\in \pm1, m\in \mathbb{N}.
\end{eqnarray}
Note that $\delta\in \Gamma^{J}$ indicates $L_{\delta}\in J$.   By Lemma 2.11, $\psi(L_{\delta})=L_{\delta\psi^{-1}}$. Since $L$ is of maximal length, we have  $\psi(L_{\delta}) \subset \psi(J)=J$. So $L_{\delta\psi^{-1}}\subset I$. Similarly we get
\begin{eqnarray}
L_{\delta\psi^{-m}}\in J, ~\mbox{for}~m\in \mathbb{N}.
\end{eqnarray}
Hence we can argue as above with the
root-multiplicativity and maximal length of $L$ from $\delta$ instead of $\delta_0$, to get that in case $\epsilon \delta_0\phi^{-m}\in \Gamma^{J}$
for some $\epsilon\in  \pm1$, then $0\neq L_{\epsilon\delta_0\psi^{-m}}\in I$.

The decomposition of $J$ in (5.12) finally gives us $I=J$.

Now suppose there is not any $\delta_0\in \Gamma^{J I}$ such that $-\delta_0\in \Gamma^{J I}$.  Then we have
  \begin{eqnarray}
\Gamma^{J}=\Gamma^{J I}\dot{\cup} -\Gamma^{J I},
\end{eqnarray}
where $-\Gamma^{J I}:=\{-\gamma|\gamma\in\Gamma^{J I}\}$.
Define
  \begin{eqnarray}
I':=(\sum_{-\gamma\in -\Gamma^{ J I},\gamma\in\Lambda}A_{\gamma}L_{-\gamma})\oplus(\bigoplus_{-\gamma\in -\Gamma^{J I}}L_{-\gamma}).
\end{eqnarray}

First, we claim  that  $I'$ is a Hom-Leibniz-ideal of $L$.
In fact, By Lemma 2.11, $\psi(L_{-\gamma})\subset L_{-\gamma\psi^{-1}}$, $-\gamma\psi^{-1}\in -\Gamma^{J I}$ and $\psi(A_{\gamma}L_{-\gamma})\subset \psi(L_0)\subset L_0$  if $A_{\gamma}L_{-\gamma}\neq 0$ (otherwise is trivial).
So  $\psi(I')\subset I'$.

Since $A_{\gamma}L_{-\gamma}\subset L_0$, by (5.15), we have
  \begin{eqnarray}
&&[L,I']\nonumber\\
&=&[H\oplus(\bigoplus_{\delta\in\Gamma^{\neg J}}L_\delta),
(\sum_{-\gamma\in -\Gamma^{ J I},\gamma\in\Lambda}A_{\gamma}L_{-\gamma})\oplus(\bigoplus_{-\gamma\in -\Gamma^{J I}}L_{-\gamma})]\subset\nonumber\\
&&[\bigoplus_{\delta\in\Gamma^{\neg J}}L_\delta,(\sum_{-\gamma\in -\Gamma^{ J I},\gamma\in\Lambda}A_{\gamma}L_{-\gamma})]
        +[\bigoplus_{\delta\in\Gamma}L_\delta,\bigoplus_{-\gamma\in -\Gamma^{J I}}L_{-\gamma}]
        +\sum_{-\gamma\in -\Gamma^{J I}}L_{-\gamma}.~~
\end{eqnarray}

For the expression $[\bigoplus_{\delta\in\Gamma^{\neg J}}L_\delta,(\sum_{-\gamma\in -\Gamma^{J I},\gamma\in\Lambda}A_{\gamma}L_{-\gamma})]$ in   (5.16),
if some $[L_\delta,A_{\gamma}L_{-\gamma}]\neq 0,$
we have that in case $\delta=-\gamma$,
$[L_{-\gamma},A_{\gamma}L_{-\gamma}]\subset  L_{-\gamma\psi^{-1}}\subset I'$,
and in case $\delta=\gamma$, since $I$ is a Hom-Leibniz-ideal of $L$, $-\gamma\notin \Gamma^{J I}$ implies
$[L_{-\gamma},A_{-\gamma}L_{\gamma}]=0$. By the  maximal length of $L$ and the symmetry of $\Gamma$,
we have $[L_{\gamma},A_{\gamma}L_{-\gamma}]=0$.
Suppose $\delta\notin \{\gamma,-\gamma\}.$
 By Definition 2.3,
 \begin{eqnarray*}
[L_\delta,A_{\gamma}L_{-\gamma}]
\subset \phi(A_{\gamma})[L_\delta, L_{-\gamma}]+\rho(L_\delta)(A_{\gamma})L_{-\gamma}.
\end{eqnarray*}
Since $(L,A)$ is regular, $\phi(A_{\gamma})\subset A_{\gamma}.$
As $[L_\delta,A_{\gamma}L_{-\gamma}]\neq 0,$  we get
$A_{\gamma}[L_\delta, L_{-\gamma}]\neq 0$ or $\rho(L_\delta)(A_{\gamma})L_{-\gamma}\neq 0$.
By the  maximal length of $L$, either $A_{\gamma}[L_\delta, L_{-\gamma}]=L_{\gamma+(\delta-\gamma)\psi^{-1}}$
or $\rho(L_\delta)(A_{\gamma})L_{-\gamma}=L_{\delta}$.
In both cases, since $\gamma\in\Gamma^{J I}$, by the root-multiplicativity of $L$,
we have $L_{-\delta}\subset I$ and therefore $-\delta\in\Gamma^{J I}$. That is, $L_{\delta}\subset I'$.
So $[\bigoplus_{\delta\in\Gamma^{\neg J}}L_\delta,(\sum_{-\gamma\in -\Gamma^{J I},\gamma\in\Lambda}A_{\gamma}$\\$L_{-\gamma})]\subset I'$.

For the expression $[\bigoplus_{\delta\in\Gamma^{\neg J}}L_\delta,\bigoplus_{-\gamma\in -\Gamma^{J I}}L_{-\gamma}]$ in (5.16),
if some $[L_\delta,L_{-\gamma}]\neq 0,$ then $[L_\delta,L_{-\gamma}]=L_{(\delta-\gamma)\psi^{-1}}$.
On the one hand, let $\delta-\gamma\neq 0$.
Since  $\gamma\in\Gamma^{J I}$, by the root-multiplicativity of $L$,
we have $[L_\gamma,L_{-\delta}]=L_{(\gamma-\delta)\psi^{-1}}\subset I$.
So $(\delta-\gamma)\psi^{-1}\in\Gamma_I$ and therefore $L_{(\delta-\gamma)\psi^{-1}}\subset I'$.
On the other hand, let $\delta-\gamma=0$. Suppose $[L_\gamma,L_{-\gamma}]\neq 0$, since $\gamma\in\Gamma^{J I}$,
we get $[L_\gamma,L_{-\gamma}]\subset I$.
Thus $L_{-\gamma}=[[L_\gamma,L_{-\gamma}],L_{-\gamma\psi}]\subset I$.
According to the discussion above, $\gamma,-\gamma\in\Gamma^{J I}$,
 a contradiction with (5.14).
So $[\bigoplus_{\delta\in\Gamma^{\neg J}}L_\delta,\bigoplus_{-\gamma\in -\Gamma^{J I}}L_{-\gamma}]\subset I'$.

Second, we claim that $\rho(I')(A)L\subset I'$.
In fact, by Definition 2.1, we have
 \begin{eqnarray*}
\rho(I')(A)L\subset [I',AL]+A[I',L]
 \end{eqnarray*}
Since $I'$ is a Hom-Leibniz-ideal of $L$, we get $[I',AL]\subset I', [I',L]\subset I'$.
So it is sufficient to verify that $AI'\subset I'$. For this, we calculate
 \begin{eqnarray}
AI'&=&(A_0\oplus(\bigoplus_{\alpha\in\Lambda}A_\alpha))
(\sum_{-\gamma\in -\Gamma^{J I},\gamma\in\Lambda}A_{\gamma}L_{-\gamma})\oplus(\bigoplus_{-\gamma\in -\Gamma^{J I}}L_{-\gamma})\subset\nonumber\\
&&I'+(\bigoplus_{\alpha\in\Lambda}A_\alpha)(\sum_{-\gamma\in -\Gamma^{J I},\gamma\in\Lambda}A_{\gamma}L_{-\gamma})
+(\bigoplus_{\alpha\in\Lambda}A_\alpha)(\bigoplus_{-\gamma\in -\Gamma^{J I}}L_{-\gamma}).
 \end{eqnarray}

 For the expression $(\bigoplus_{\alpha\in\Lambda}A_\alpha)(\sum_{-\gamma\in -\Gamma^{J I},\gamma\in\Lambda}A_{\gamma}L_{-\gamma})$ in (5.17),
 if some $A_\alpha(A_{\gamma}L_{-\gamma})\neq 0$, we have that in case $\alpha=-\gamma$, clearly
 $A_\alpha(A_{\gamma}L_{-\gamma})=A_{-\gamma}(A_{\gamma}L_{-\gamma})\subset L_{-\gamma}\subset I'$.
 In case of  $\alpha=\gamma$, since $-\gamma\notin \Gamma_I$, we get $A_{-\gamma}(A_{-\gamma}L_{\gamma})=0$.
 By the   the  maximal length of $L$, we have $A_\alpha(A_{\gamma}L_{-\gamma})=A_{\gamma}(A_{\gamma}L_{-\gamma})=0$.
 Suppose that $\alpha\notin\{\gamma,-\gamma\}$,
by the   the  maximal length of $L$, we have $A_\alpha(A_{\gamma}L_{-\gamma})=L_\alpha$.
Since $\gamma\in\Gamma^{J I}$, by the root-multiplicativity of $L$, we have $L_{-\gamma}\subset I$, that is, $-\alpha\in\Gamma^{J I}$.
 So $\alpha\in-\Gamma^{J I}$ and $L_\alpha\subset I'$.
 Thus $(\bigoplus_{\alpha\in\Lambda}A_\alpha)(\sum_{-\gamma\in -\Gamma^{J I},\gamma\in\Lambda}A_{\gamma}L_{-\gamma})\subset I'$.

  For the expression $(\bigoplus_{\alpha\in\Lambda}A_\alpha)(\bigoplus_{-\gamma\in -\Gamma^{J I}}L_{-\gamma})$ in (5.17),
if some $A_\alpha L_{-\gamma}\neq 0$, in case $\alpha-\gamma\in \Gamma^{J I}$,
by the root-multiplicativity of $L$, we have $A_{-\alpha} L_{\gamma}\neq 0$.
Again by the  maximal length of $L$, we have $A_{-\alpha} L_{\gamma}=L_{-\alpha+\gamma}$.
So $ -\alpha+\gamma\in \Gamma^{J I}$, a contradiction. Thus $ \alpha-\gamma\in -\Gamma^{J I}$ and therefore
$(\bigoplus_{\alpha\in\Lambda}A_\alpha)(\bigoplus_{-\gamma\in -\Gamma^{J I}}L_{-\gamma})\subset I'$.

By the discussion above, we have shown that $\rho(I')(A)L\subset I'$ and therefore $I'$ is an ideal of $(L,A)$.

Finally, we will verify that $L=I\oplus I'$ with ideals $I,I'$.  Since $[I',I]=0$, by the commutativity of  $H$, we get $\sum_{\gamma\in\Gamma}[L_{\gamma},L_{-\gamma}]=0$,
so $H$ must has the form
  \begin{eqnarray}
H=(\sum_{\gamma\in \Gamma^{J I},-\gamma\in\Lambda}A_{\gamma}L_{-\gamma})\oplus(\sum_{-\gamma\in -\Gamma^{J I},\gamma\in\Lambda}A_{\gamma}L_{-\gamma}).
\end{eqnarray}
In order to show that the sum in (5.19) is direct, we take any
$h\in(\sum_{\gamma\in \Gamma^{J I},-\gamma\in\Lambda}A_{\gamma}L_{-\gamma})\oplus(\sum_{-\gamma\in -\Gamma^{J I},\gamma\in\Lambda}A_{\gamma}L_{-\gamma}).$
Suppose $h\neq 0$, then $h\notin Z(L)$.
Since $L$ is split, there is $v_\delta\in L_\delta, \delta\in \Gamma$ satisfying
$[h,v_\delta]=\delta(h)\psi(v_\delta)\neq 0$.
By Proposition 3.3, $0\neq\delta(h)\psi(v_\delta)\in L_{\delta\psi^{-1}}$.
While $L_{\delta\psi^{-1}}\subset I\cap I'=0$, a contradiction.
So $h=0$, as required. And this finishes the proof.
$\hfill \Box$

\begin{corollary}
Let $(L,A)$ be a tight split regular Hom-Leibniz-Rinehart algebra of  maximal length and root-multiplicative.
If $\Gamma^{J}$, $\Gamma^{\neg J}$ are symmetric and $\Gamma^{\neg J}$ has all of its roots $\neg J$-connected, and $\Lambda$ have all its nonzero weights connected.
Then
\begin{eqnarray*}
L=\bigoplus_{i\in I}L_i,~~A=\bigoplus_{j\in K}A_j,
\end{eqnarray*}
where any $L_i$ is a simple ideal of $L$ having all of its nonzero roots connected satisfying $[L_i,L_{i'}]=0$ for any $i'\in I$ with $i\neq i'$,
and any $A_j$ is a simple ideal of $A$ satisfying $A_j A_{j'}=0$ for any $j'\in K$ such that $j'\neq j$.

Furthermore, for any $i\in I$ there exists a unique $\overline{j}\in K$ such that $A_{\overline{j}}L_i\neq 0$.
We also have that any $L_i$ is a split regular Hom-Leibniz-Rinehart algebra over $A_{\overline{j}}$.
\end{corollary}

{\bf Proof.}
  It is analogous to Theorem 5.7 in \cite{Wang19}.
$\hfill \Box$

 \begin{center}
 {\bf ACKNOWLEDGEMENT}
 \end{center}

 The paper is supported by  the NSF of China (Nos. 11761017 and 11801304),
 the Youth Project for Natural Science Foundation of Guizhou provincial department of education (No. KY[2018]155)
 and the Anhui Provincial Natural Science Foundation (Nos. 1908085MA03 and 1808085MA14).

\end{document}